\documentclass[12pt,aps,onecolumn,nobibnotes]{revtex4}

\usepackage{graphicx}%
\usepackage{dcolumn}
\usepackage{amsmath}
\usepackage{amssymb}

\makeatletter
\def\btt#1{\texttt{\@backslashchar#1}}%
\DeclareRobustCommand\bblash{\btt{\@backslashchar}}%
\makeatother
\def\thepart          {\Roman{part}} %
\def\thesection       {\arabic{section}}
\def\thesubsection    {\thesection.\arabic{subsection}.}
\def\theparagraph     {\alph{paragraph}}

\def\cop{\Delta}

\def\cnt{\varepsilon}
\def\ot{\otimes}

\def\GLqtwo{{GL}_{q}(2)}
\def\GLhtwo{{GL}_{h}(2)}

\def\GLpqtwo{{GL}_{p,q}(2)}
\def\GLhh'two{{GL}_{h,h'}(2)}

\def\Grss'{{G}_{r}^{s,s'}}
\def\Gmkk'{{G}_{m}^{k,k'}}

\def\ident{{\bf 1}}

\def\ca{{\mathcal A}}

\def\ch{{\mathcal H}}

\def\cu{{\mathcal U}}

\def\ms{\mathsf s}

\begin{document}

\title[Short Title]{COLOURED EXTENSION OF $GL_{q}(2)$ AND ITS DUAL
ALGEBRA}
\author{Deepak Parashar}
\email{Deepak.Parashar@mis.mpg.de}
\affiliation{
Max-Planck-Institute for Mathematics in the Sciences \\
Inselstrasse 22-26, D-04103 Leipzig \\
Germany
}
\begin{abstract}
We address the problem of duality between the coloured extension of the
quantised algebra of functions on a group and that of its quantised
universal enveloping algebra {\it i.e.} its dual. In particular, we derive
explicitly the algebra dual to the coloured extension of $GL_{q}(2)$ using
the coloured $RLL$ relations and exhibit its Hopf structure. This leads to
a coloured generalisation of the $R$-matrix procedure to construct a
bicovariant differential calculus on the coloured version of
$GL_{q}(2)$. In addition, we also propose a coloured generalisation of the
geometric approach to quantum group duality given by Sudbery and Dobrev.
\end{abstract}

\maketitle

\section{Introduction}

The quantum group $\GLqtwo$ is known to admit a coloured extension by
introducing some continuously varying colour parameters associated to the
generators. In such an extension, the associated algebra and the coalgebra
are defined in a way that all Hopf algebraic properties remain
preserved. Such extensions have been introduced in \cite{ge,bh,oht} and
studied by various authors \cite{kb,jagan,quesne,preeti} in recent 
years. However, some of the basic algebro-geometric structure underlying
these coloured extensions still needs to be established. As
such, we shall focus on the coloured extension of the most intuitive 
quantum group $\GLqtwo$. While some aspects of this example have already
been studied from both, the standard $q$-deformations as well as the
Jordanian (nonstandard) $h$-deformations \cite{quesne,preeti}, it has only
recently been shown \cite{deeps} that the contraction procedure could be
used to obtain the coloured Jordanian quantum groups from their coloured
$q$-deformed counterparts. In particular, the coloured extension of
$\GLqtwo$ was treated in detail in \cite{deeps} to obtain a new coloured
extension of Jordanian $\GLhtwo$.
\\
In the present paper, we investigate the algebra dual to the coloured
extension of $\GLqtwo$ by generalising two well-known approaches to the
problem: the (algebraic) $R$-matrix approach \cite{frt} and the geometric
approach \cite{sud,dob} due to Sudbery and Dobrev. We first clarify the
notion of duality between a coloured quantum group and its dual {\em i.e.}
the coloured quantised universal enveloping algebra. We then generalise
the $R$-matrix approach to establish duality for the coloured extension of
$\GLqtwo$ and we obtain a new coloured quantum algebra corresponding to
$gl(2)$ and exhibit its Hopf algebra structure. The coloured $R$-matrix
procedure naturally leads us to formulate a constructive differential
calculus \cite{jurco} on the coloured extension of $\GLqtwo$.
\\
Furthermore, we propose a coloured generalisation of the geometric notion
of duality for quantum groups {\em i.e.} regarding the dual algebra as the
algebra of tangent vectors at the identity of the group. This
generalisation could also be of significance in establishing the duality
for the coloured extension of Jordanian quantum groups.

\section{Coloured extension of $GL_{{q}}(2)$}

The coloured extension of the quantum group $\GLqtwo$ is governed by the
coloured $R$-matrix \cite{kb},
\begin{equation}
R_{q}^{\lambda,\mu} = \begin{pmatrix}
q^{1-(\lambda-\mu)} & 0 & 0 & 0\\
0 & q^{\lambda+\mu} & 0 & 0\\
0 & q-q^{-1} & q^{-(\lambda+\mu)} & 0\\
0 & 0 & 0 & q^{1+(\lambda-\mu)}\
\end{pmatrix}
\end{equation}
which  is nonadditive {\em i.e.} $R^{\lambda,\mu} \neq R(\lambda - \mu)$.
It satisfies the so-called coloured quantum Yang-Baxter equation
\begin{equation}
R_{12}^{\lambda,\mu}R_{13}^{\lambda,\nu}R_{23}^{\mu,\nu} =
R_{23}^{\mu,\nu}R_{13}^{\lambda,\nu}R_{12}^{\lambda,\mu}
\end{equation}
which is in general
multicomponent and $\lambda$, $\mu$, $\nu$ are considered
as `colour' variables. The $RTT$ relations are also extended to
incorporate the coloured extension as
\begin{equation}
R_{q}^{\lambda,\mu}T_{1\lambda}T_{2\mu}=T_{2\mu}T_{1\lambda}
R_{q}^{\lambda,\mu}
\end{equation}
(where $T_{1\lambda}=T_{\lambda}\ot \ident$ and
$T_{2\mu}=\ident\ot T_{\mu}$) in which the entries of the $T$
matrices carry colour  dependence {\em i.e.}
$T_{\lambda}= \left(
\begin{smallmatrix}a_{\lambda}&b_{\lambda}\\c_{\lambda}&d_{\lambda}
\end{smallmatrix}\right)$,
$T_{\mu}= \left(
\begin{smallmatrix}a_{\mu}&b_{\mu}\\c_{\mu}&d_{\mu}\end{smallmatrix} 
\right)$.
The coproduct and counit for the coalgebra structure are given by
$\cop (T_{\lambda})=T_{\lambda}\dot{\ot} T_{\lambda}$,
$\cnt (T_{\lambda})=\ident$.
The quantum determinant
$D_{\lambda}=a_{\lambda}d_{\lambda}-r^{-(1+2\lambda)}c_{\lambda}b_{\lambda}$
is group-like but not central. The antipode is given by
\begin{equation}
S(T_{\lambda})= 
D_{\lambda}^{-1}\begin{pmatrix}d_{\lambda}&-r^{1+2\lambda}b_{\lambda}\\
-r^{-1-2\lambda}c_{\lambda}&a_{\lambda}\end{pmatrix}
\end{equation}
and depends on one colour variable at a time. The full Hopf
algebraic structure can be constructed resulting in a coloured extension
of $\GLqtwo$ within the framework of the $FRT$ formalism. Since $\lambda$
and $\mu$ are continuous variables, this implies the coloured extension of
$\GLqtwo$ has an infinite number of generators. The colourless limit
$\lambda = \mu = 0$ gives back the ordinary single-parameter deformed
quantum group $\GLqtwo$, and the monochromatic limit $\lambda = \mu \neq
0$ gives rise to the uncoloured two-parameter deformed quantum group
$\GLpqtwo$.

\section{Duality (R-matrix approach)}

In this section, we investigate in detail the dual structure for the
coloured extension of $\GLqtwo$ employing the $R$-matrix approach. In
doing so, let us denote the generators of the yet unknown dual
algebra by 
$\{ A_{\lambda},B_{\lambda},C_{\lambda},D_{\lambda} \}$ and
$\{ A_{\mu},B_{\mu},C_{\mu},D_{\mu} \}$. The following
pairings hold
\begin{equation}
\langle A_{\lambda | \mu}, a_{\lambda | \mu} \rangle = 
\langle B_{\lambda | \mu}, b_{\lambda | \mu} \rangle =
\langle C_{\lambda | \mu}, c_{\lambda | \mu} \rangle =
\langle D_{\lambda | \mu}, d_{\lambda | \mu} \rangle = 
\ident
\end{equation}
All other pairings give zeroes and the notation $\lambda|\mu$ in the
subscript in the above relations means {\em either} $\lambda$ {\em or}
$\mu$. The $R^{+}$ and $R^{-}$ matrices corresponding to the coloured extension
of $\GLqtwo$ are
\begin{eqnarray}  
R^{+} &=& c^{+} q^{1/2} \begin{pmatrix}
q^{-1/2}q^{1-\lambda+\mu} & 0 & 0 & 0\\
0 & q^{-1/2}q^{-(\lambda+\mu)} & q^{-1/2}(q-q^{-1}) & 0\\
0 & 0 & q^{-1/2}q^{\lambda+\mu} & 0\\
0 & 0 & 0 & q^{-1/2}q^{1+\lambda-\mu} 
\end{pmatrix} \\
R^{-} &=& c^{-} q^{-1/2} \begin{pmatrix}
q^{1/2}q^{-(1-\lambda+\mu)} & 0 & 0 & 0\\
0 & q^{1/2}q^{-(\lambda+\mu)} & 0 & 0\\  
0 & -q^{1/2}(q-q^{-1}) & q^{1/2}q^{\lambda+\mu} & 0\\
0 & 0 & 0 & q^{1/2}q^{-(1+\lambda-\mu)}
\end{pmatrix}
\end{eqnarray} 
where $R^{+}=c^{+}R_{21}$ and $R^{-}=c^{-}R_{12}^{-1}$ by definition.
The coloured $L^{\pm}$ functionals can be expressed as
\begin{eqnarray}
L^{+}_{\lambda(\mu)} &=& c^{+} q^{1/2} \begin{pmatrix}
q^{H_{\lambda(\mu)}/2}q^{\mu H_{\lambda(\mu)}-\lambda H'_{\lambda(\mu)}} &
q^{-1/2}(q-q^{-1})C_{\lambda(\mu)}\\
0 & q^{-H_{\lambda(\mu)}/2}q^{\mu H_{\lambda(\mu)}+\lambda
H'_{\lambda(\mu)}}
\end{pmatrix}\\
L^{-}_{\lambda(\mu)} &=& c^{-} q^{-1/2} \begin{pmatrix}
q^{-H_{\lambda(\mu)}/2}q^{\lambda H_{\lambda(\mu)}-\mu H'_{\lambda(\mu)}}
& 0\\
q^{1/2}(q^{-1}-q)B_{\lambda(\mu)} & 
q^{H_{\lambda(\mu)}/2}q^{\lambda H_{\lambda(\mu)}+\mu H'_{\lambda(\mu)}}
\end{pmatrix}
\end{eqnarray}
where $H_{\lambda}=A_{\lambda}-D_{\lambda}$,
$H'_{\lambda}=A_{\lambda}+D_{\lambda}$ and
$H_{\mu}=A_{\mu}-D_{\mu}$, $H'_{\mu}=A_{\mu}+D_{\mu}$.
The notation $\lambda(\mu)$ in the subscript means $\lambda$
({\em respectively} $\mu$). So, $L^{+}_{\lambda(\mu)}$ means
$L^{+}_{\lambda}$ ({\em resp.} $L^{+}_{\mu}$)
Each one of $L^{\pm}_{\lambda}$ and $L^{\pm}_{\mu}$ depends on
both $\lambda$ and $\mu$. The notation $L^{\pm}_{\lambda}$
implies that the generators of the dual carry $\lambda$ dependence, and
similarly $L^{\pm}_{\mu}$ implies that the generators of the dual carry
$\mu$ dependence. The duality pairings are then given by the action of
the functionals $L^{\pm}_{\lambda}$ and $L^{\pm}_{\mu}$ on the
$T$-matrices $T_{\lambda}$ and $T_{\mu}$
\begin{eqnarray}
&&(L^{+}_{\lambda|\mu})^{a}_{b} (T_{\lambda|\mu})^{c}_{d} =
(R^{+})^{ac}_{bd}\\
&&(L^{-}_{\lambda|\mu})^{a}_{b} (T_{\lambda|\mu})^{c}_{d} =
(R^{-})^{ac}_{bd}
\end{eqnarray}
Again, according to the notation introduced $T_{\lambda|\mu}$ implies
$T_{\lambda}$ {\em or} $T_{\mu}$ and $L^{\pm}_{\lambda|\mu}$ implies
$L^{\pm}_{\lambda}$ {\em or} $L^{\pm}_{\mu}$. For vanishing colour
variables, the coloured $L^{\pm}$ functionals reduce to the ordinary
$L^{\pm}$ functionals for $\GLqtwo$. The
commutation relations of the algebra dual to a coloured quantum group can
be obtained from the modified or the {\em coloured} $RLL$ relations
\begin{eqnarray}
R_{12} L^{\pm}_{2\lambda} L^{\pm}_{1\mu} &=& 
L^{\pm}_{1\mu} L^{\pm}_{2\lambda} R_{12}\\
R_{12} L^{+}_{2\lambda} L^{-}_{1\mu} &=&
L^{-}_{1\mu} L^{+}_{2\lambda} R_{12}
\end{eqnarray}
using the coloured $L^{\pm}$ functionals where
$L^{\pm}_{1\mu}=L^{\pm}_{\mu}\ot \ident$ and
$L^{\pm}_{2\lambda}=\ident \ot L^{\pm}_{\lambda}$. Using the above
formulae, we obtain the commutation relations between the generating
elements of the algebra dual to the coloured extension of $\GLqtwo$.
\begin{equation}
\begin{array}{lll}
\left[A_{\lambda},B_{\mu}\right]=B_{\mu} \qquad &
\left[D_{\lambda},B_{\mu}\right]=-B_{\mu} \qquad &\\
\left[A_{\lambda},C_{\mu}\right]=-C_{\mu} \qquad &
\left[D_{\lambda},C_{\mu}\right]=C_{\mu} \qquad &\\
\left[A_{\lambda},D_{\mu}\right]=0 \qquad &
\left[H_{\lambda},H_{\mu}\right]=0 \qquad &
\left[H'_{\lambda},\bullet\right]=0
\end{array}
\end{equation}
\begin{equation}
q^{-(\lambda + \mu)}C_{\lambda}B_{\mu} - 
q^{\lambda + \mu}B_{\mu}C_{\lambda} = 
\frac{q^{\lambda H_{\mu} + \mu H_{\lambda}}}{q-q^{-1}}
\left[
q^{-\frac{1}{2}(H_{\lambda}+H_{\mu})} 
q^{\lambda H'_{\lambda} - \mu H'_{\mu}} -
q^{\frac{1}{2}(H_{\lambda}+H_{\mu})}
q^{-\lambda H'_{\lambda} + \mu H'_{\mu}}
\right]
\end{equation}
\begin{equation}
\begin{array}{lll}
A_{\lambda}A_{\mu} &=& A_{\mu}A_{\lambda}\\
B_{\lambda}B_{\mu} &=& q^{2(\mu-\lambda)}B_{\mu}B_{\lambda}\\
C_{\lambda}C_{\mu} &=& q^{2(\lambda-\mu)}C_{\mu}C_{\lambda}\\
D_{\lambda}D_{\mu} &=& D_{\mu}D_{\lambda}
\end{array}
\end{equation}
where $H_{\lambda}$ and $H'_{\lambda}$ are as before. The relations
satisfy the $\lambda \leftrightarrow \mu$ exchange symmetry. The
associated coproduct of the elements of the dual algebra is given by
\begin{eqnarray}
\cop(A_{\lambda(\mu)}) &=& A_{\lambda(\mu)}\ot \ident + \ident\ot
A_{\lambda(\mu)}\\
\cop(B_{\lambda(\mu)}) &=& B_{\lambda(\mu)}\ot
q^{A_{\lambda(\mu)}-D_{\lambda(\mu)}} + \ident\ot B_{\lambda(\mu)}\\
\cop(C_{\lambda(\mu)}) &=& C_{\lambda(\mu)}\ot
q^{A_{\lambda(\mu)}-D_{\lambda(\mu)}} + \ident\ot C_{\lambda(\mu)}\\
\cop(D_{\lambda(\mu)}) &=& D_{\lambda(\mu)}\ot \ident + \ident\ot
D_{\lambda(\mu)}
\end{eqnarray}
The counit $\cnt(Y_{\lambda|\mu}) = 0;$ where $Y_{\lambda(\mu)} = \{
A_{\lambda(\mu)}, B_{\lambda(\mu)}, C_{\lambda(\mu)}, D_{\lambda(\mu)} \}$
and the antipode is
\begin{eqnarray}
S(A_{\lambda(\mu)}) &=& - A_{\lambda(\mu)}\\
S(B_{\lambda(\mu)}) &=& -
B_{\lambda(\mu)}q^{-(A_{\lambda(\mu)}-D_{\lambda(\mu)})}\\
S(C_{\lambda(\mu)}) &=& -
C_{\lambda(\mu)}q^{-(A_{\lambda(\mu)}-D_{\lambda(\mu)})}\\
S(D_{\lambda(\mu)}) &=& - D_{\lambda(\mu)}
\end{eqnarray}
Thus we have defined a new single-parameter coloured quantum algebra
corresponding to $gl(2)$, which in the monochromatic limit defines the
standard uncoloured two-parameter quantum algebra for $gl(2)$.

\section{Constructive calculus}

We now proceed towards a coloured generalisation of the constructive
differential calculus \cite{jurco,cast} for the coloured extension of
$GL_{q}(2)$. Analogous to the standard uncoloured
quantum group, a bimodule $\Gamma$ (space of quantum one-forms
$\omega$) is characterised by the commutation relations between $\omega$
and $a_{\lambda(\mu)} \in\ca$, the {\em coloured} quantum group
corresponding to $GL_{q}(2)$
\begin{equation}
\omega a_{\lambda(\mu)} = (\ident\ot f_{\lambda,\mu}) \cop
(a_{\lambda(\mu)}) \omega
\end{equation}
and the linear functional $f_{\lambda,\mu}$ is defined in terms of the
coloured $L^{\pm}$ matrices
\begin{equation}
f_{\lambda,\mu} = S(L_{\lambda|\mu}^{+})L_{\lambda|\mu}^{-}
\end{equation}
Thus, we have
\begin{equation}
\omega a_{\lambda(\mu)} = [(\ident\ot
S(L_{\lambda|\mu}^{+})L_{\lambda|\mu}^{-}) \cop (a_{\lambda(\mu)})]
\omega
\end{equation}
In terms of components, this can be written as
\begin{equation}
\omega_{ij} a_{\lambda(\mu)} = [(\ident\ot
S(l^{+}_{(\lambda|\mu)ki})l^{-}_{(\lambda|\mu)jl}) \cop
(a_{\lambda(\mu)})] \omega_{kl}
\end{equation}
using the expressions $L^{\pm} = l^{\pm}_{ij}$ and $\omega =
\omega_{ij}$ where $i, j = 1,2$. From these relations, one can obtain the
commutation relations of all the left-invariant one-forms with the
elements of the coloured extension of $GL_{q}(2)$
The left-invariant vector fields $\chi_{ij}$ on $\ca$ are given by the
expression
\begin{equation}
\chi_{ij} = S(l^{+}_{(\lambda|\mu)ik})l^{-}_{(\lambda|\mu)kj} -
\delta_{ij}\cnt
\end{equation}
The vector fields act on the elements $a_{\lambda(\mu)}$ of the coloured 
quantum group as
\begin{equation}
\chi_{ij} a_{\lambda(\mu)} =
(S(l^{+}_{(\lambda|\mu)ik})l^{-}_{(\lambda|\mu)kj} -
\delta_{ij}\cnt)a_{\lambda(\mu)}
\end{equation}
Furthermore, using the formula
$\mathbf{d} a_{\lambda(\mu)} =\sum_{i}(\chi_{i} \ast
a_{\lambda(\mu)}) \omega^{i}$, we obtain the action of the exterior
derivative on the generating elements
\begin{eqnarray}
\mathbf{d} a_{\lambda(\mu)} &=&
(\ms q^{-2+2(\lambda-\mu)}-1)a_{\lambda(\mu)}\omega^{1} +
\ms (q^{-1}-q)q^{\lambda+\mu}b_{\lambda(\mu)}\omega^{+}\\
& &+ (\ms -1)a_{\lambda(\mu)}\omega^{2} \notag \\
\mathbf{d} b_{\lambda(\mu)} &=& 
(\ms (q^{-1}-q)^{2}+\ms -1)b_{\lambda(\mu)}\omega^{1} +
\ms (q^{-1}-q)q^{-(\lambda+\mu)}a_{\lambda(\mu)}\omega^{-}\\
& &+ (\ms q^{-2+2(\mu-\lambda)}-1)b_{\lambda(\mu)}\omega^{2} \notag \\
\mathbf{d} c_{\lambda(\mu)} &=&
(\ms q^{-2+2(\lambda-\mu)}-1)c_{\lambda(\mu)}\omega^{1} +
\ms (q^{-1}-q)q^{\lambda+\mu}d_{\lambda(\mu)}\omega^{+}\\
& &+ (\ms -1)c_{\lambda(\mu)}\omega^{2} \notag \\
\mathbf{d} d_{\lambda(\mu)} &=&
(\ms (q^{-1}-q)^{2}+\ms -1)d_{\lambda(\mu)}\omega^{1} +
\ms (q^{-1}-q)q^{-(\lambda+\mu)}c_{\lambda(\mu)}\omega^{-}\\ 
& &+ (\ms q^{-2+2(\mu-\lambda)}-1)d_{\lambda(\mu)}\omega^{2} \notag
\end{eqnarray}
where $\omega^{1}=\omega_{11}$, $\omega^{+}=\omega_{12}$,
$\omega^{-}=\omega_{21}$, $\omega^{2}=\omega_{22}$ and
$\ms=(c^{+})^{-1}c^{-}$. $\mathbf{d}\ca$ generates $\Gamma$ as a left
$\ca$-module. This defines a first order differential calculus
$(\Gamma, \mathbf{d})$ on the coloured extension of $\GLqtwo$. Since
the colour variables $\lambda$ and $\mu$ are continuously varying, the
differential calculus obtained is infinite-dimensional. The differential
calculus on the uncoloured single-parameter quantum group $\GLqtwo$ is
easily recovered in the colourless limit and that of the uncoloured
two-parameter quantum group $\GLpqtwo$ in the monochromatic limit.

\section{Dual basis (Geometric Approach)}

It is well-known that two bialgebras $\cu$ and $\ca$ are in duality if
there exists a doubly nondegenerate bilinear form
\begin{equation}
\langle,\rangle:\cu \ot \ca\rightarrow \mathbf{C};\quad
\langle,\rangle:(u,a)\rightarrow \langle u,a\rangle;\quad
\forall u\in \cu, a\in \ca
\end{equation}
such that for $u,v\in \cu$ and $a,b\in \ca$, we have
\begin{equation}
\begin{array}{ll}
\langle u,ab\rangle = \langle \cop_{\cu}(u), a\ot b\rangle\\
\langle uv,a\rangle = \langle u\ot v, \cop_{\ca}(a) \rangle
\end{array}
\end{equation}
\begin{equation}
\begin{array}{ll}
\langle \ident_{\cu},a\rangle=\cnt_{\ca}(a)\\
\langle u,\ident_{\ca}\rangle=\cnt_{\cu}(u)
\end{array}
\end{equation}
For the two bialgebras to be in duality as Hopf algebras, $\cu$ and
$\ca$ further satisfy
\begin{equation}
\langle S_{\cu}(u),a \rangle=\langle u,S_{\ca}(a)\rangle
\end{equation}
It is enough to define the pairing between the generating elements of
the two algebras. Pairing for any other elements of $\cu$ and $\ca$
follows from these relations and the bilinear form inherited by the
tensor product. The geometric approach for duality for quantum groups
was motivated by the fact that, at the classical level, an element of the
Lie algebra corresponding to a Lie group is a tangent vector at the
identity of the Lie group. Let $\ch$ be a given Hopf algebra generated by
non-commuting elements $a$, $b$, $c$, $d$. The $q$-analogue of tangent
vector at the identity would then be obtained by first differentiating the
elements of the given Hopf algebra $\ch$ (polynomials in $a$, $b$, $c$,
$d$) and then putting  
$\left( \begin{smallmatrix}a&b\\c&d\end{smallmatrix} \right) =
\left( \begin{smallmatrix}1&0\\0&1\end{smallmatrix} \right)$
later on (i.e., taking the counit operation analogous to the unit
element at the group level). The elements thus obtained
would belong to the dual Hopf algebra $\ch^{*}$. The approach is due to
Sudbery \cite{sud} and Dobrev \cite{dob} and has proved to be quite a
powerful tool in understanding the quantum group duality from a geometric
point of view. In what follows in this section, we propose to give a
coloured generalisation of such a geometric picture of duality using the
example of $GL_{q}(2)$. Let $\ca_{q}^{\lambda,\mu}$ denote the coloured
extension  of $GL_{q}(2)$. Then, as a Hopf algebra
$\ca_{q}^{\lambda,\mu}$ is generated by elements $y_{\lambda}=\{
a_{\lambda},b_{\lambda},c_{\lambda},d_{\lambda} \}$ and
$y_{\mu}=\{ a_{\mu},b_{\mu},c_{\mu},d_{\mu} \}$. The dual basis is given
by all monomials of the form
\begin{equation} 
\begin{array}{ll}
g_{\lambda} = g_{\lambda;klmn} =
a_{\lambda}^{k}d_{\lambda}^{l}b_{\lambda}^{m}c_{\lambda}^{n} \\
g_{\mu} = g_{\mu;klmn} =
a_{\mu}^{k}d_{\mu}^{l}b_{\mu}^{m}c_{\mu}^{n}
\end{array}
\end{equation}
where  
$k,l,m,n\in Z_{+}$, and $\delta_{0000}$ is the unit of the algebra 
$\ident_{\ca}$. We use a normal ordering as follows; first put the
diagonal elements from the $T_{\lambda(\mu)}$- matrix then use the
lexicographic order
for the others. Let $\cu_{q}^{\lambda,\mu}$ be the algebra generated by
tangent vectors at the identity of $\ca_{q}^{\lambda,\mu}$. Then
$\cu_{q}^{\lambda,\mu}$ is dually paired with
$\ca_{q}^{\lambda,\mu}$. The pairing is defined through the
coloured $q$-tangent vectors as follows
\begin{eqnarray}
\langle Y_{\lambda},g_{\lambda}\rangle &=& \frac{\partial
g_{\lambda}}{\partial y_{\lambda}}|_{
  \left(      
\begin{smallmatrix} 
a_{\lambda}&b_{\lambda}\\c_{\lambda}&d_{\lambda}\end{smallmatrix}\right) =
  \left( \begin{smallmatrix}1&0\\0&1\end{smallmatrix}\right)
} =
\cnt(\frac{\partial g_{\lambda}}{\partial y_{\lambda}})
\\
\langle Y_{\mu},g_{\lambda}\rangle &=& \frac{\partial
g_{\lambda}}{\partial y_{\lambda}}|_{
  \left(
\begin{smallmatrix}
a_{\lambda}&b_{\lambda}\\c_{\lambda}&d_{\lambda}\end{smallmatrix}\right) =
  \left( \begin{smallmatrix}1&0\\0&1\end{smallmatrix}\right)  
} =
\cnt(\frac{\partial g_{\lambda}}{\partial y_{\lambda}})
\\
\langle Y_{\lambda},g_{\mu}\rangle &=& \frac{\partial
g_{\mu}}{\partial y_{\mu}}|_{
  \left(   
\begin{smallmatrix}
a_{\mu}&b_{\mu}\\c_{\mu}&d_{\mu}\end{smallmatrix}\right) =
  \left( \begin{smallmatrix}1&0\\0&1\end{smallmatrix}\right)
} =
\cnt(\frac{\partial g_{\mu}}{\partial y_{\mu}})
\\
\langle Y_{\mu},g_{\mu}\rangle &=& \frac{\partial
g_{\mu}}{\partial y_{\mu}}|_{
  \left(   
\begin{smallmatrix}
a_{\mu}&b_{\mu}\\c_{\mu}&d_{\mu}\end{smallmatrix}\right) =
  \left( \begin{smallmatrix}1&0\\0&1\end{smallmatrix}\right)
} =
\cnt(\frac{\partial g_{\mu}}{\partial y_{\mu}})
\end{eqnarray}
where $Y_{\lambda} = \left\{
A_{\lambda},B_{\lambda},C_{\lambda},D_{\lambda}\right\}$ and 
$Y_{\mu} = \left\{ 
A_{\mu},B_{\mu},C_{\mu},D_{\mu}\right\}$ 
are the sets of generating elements of the dual algebra (which has unit
$\ident_{\cu}$). More compactly, one can write
\begin{equation}
\langle Y_{\lambda | \mu},g_{\lambda(\mu)}\rangle = \cnt(\frac{\partial
g_{\lambda(\mu)}}{\partial y_{\lambda(\mu)}})
\end{equation}
Explicitly, we obtain
\begin{eqnarray}
\langle A_{\lambda|\mu},g_{\lambda(\mu)}\rangle &=& \cnt(\frac{\partial
g_{\lambda(\mu)}}{\partial a_{\lambda(\mu)}}) 
  = k\delta_{m0}\delta_{n0}\\
\langle B_{\lambda|\mu},g_{\lambda(\mu)}\rangle &=& \cnt(\frac{\partial
g_{\lambda(\mu)}}{\partial b_{\lambda(\mu)}})
   = \delta_{m1}\delta_{n0}\\
\langle C_{\lambda|\mu},g_{\lambda(\mu)}\rangle &=& \cnt(\frac{\partial
g_{\lambda(\mu)}}{\partial c_{\lambda(\mu)}})
   = \delta_{m0}\delta_{n1}\\
\langle D_{\lambda|\mu},g_{\lambda(\mu)}\rangle &=& \cnt(\frac{\partial
g_{\lambda(\mu)}}{\partial d_{\lambda(\mu)}})
   = l\delta_{m0}\delta_{n0}
\end{eqnarray}
where differentiation is from the right. As a consequence of the above
pairings, the following relations hold
\begin{eqnarray}
\langle A_{\lambda|\mu}, T_{\lambda|\mu} \rangle
&=& \left( \begin{smallmatrix} 1&0\\0&0 \end{smallmatrix}\right) \\
\langle B_{\lambda|\mu}, T_{\lambda|\mu} \rangle
&=& \left( \begin{smallmatrix} 0&1\\0&0 \end{smallmatrix}\right) \\
\langle C_{\lambda|\mu}, T_{\lambda|\mu} \rangle
&=& \left( \begin{smallmatrix} 0&0\\1&0 \end{smallmatrix}\right) \\
\langle D_{\lambda|\mu}, T_{\lambda|\mu} \rangle
&=& \left( \begin{smallmatrix} 0&0\\0&1 \end{smallmatrix}\right)
\end{eqnarray}
where
$T_{\lambda}= \left(
\begin{smallmatrix}a_{\lambda}&b_{\lambda}\\c_{\lambda}&d_{\lambda}
\end{smallmatrix}\right)$ and
$T_{\mu}= \left(
\begin{smallmatrix}a_{\mu}&b_{\mu}\\c_{\mu}&d_{\mu}\end{smallmatrix} 
\right)$ as before. Furthermore,
\begin{equation}
\langle Y_{\lambda|\mu}, \ident_{\ca}\rangle = 0; \quad 
\langle\ident_{\cu},g_{\lambda|\mu}\rangle =
\cnt_{\ca}(g_{\lambda|\mu})= \delta_{m0}\delta_{n0}
\end{equation}
The action of the monomials in $\cu_{q}^{\lambda,\mu}$ on $g_{\lambda}$
and $g_{\mu}$ then lead to the coloured $q$-commutation relations between
the generators of the dual algebra.

\section{Concluding remarks}

We have investigated the structure of the coloured extension of the
quantum group $\GLqtwo$ and its dual algebra. After establishing the
notion of duality, the dual algebra has been derived explicitly using the
$R$-matrix approach. We not only obtain a new coloured quantum algebra
corresponding to $gl(2)$ but also show that such a coloured generalisation
of the $R$-matrix approach leads to formulating a constructive
differential calculus for the coloured case. The colourless and the
monochromatic limits of both the dual algebra as well as the differential
calculus are in agreement with already known results for $\GLqtwo$ and
$\GLpqtwo$.\\
In the later part of the paper, we have proposed a generalisation of the
geometric picture of duality to incorporate the coloured extensions of
quantum groups. It would be interesting to investigate in detail this
setting in the context of the coloured Jordanian quantum groups. The
results easily extend to the higher-dimensional and multi-parametric
cases.

\end{document}